\documentclass[12pt,a4paper,twoside]{article}

\newcommand{\pageformat}[6]{\setlength{\hoffset}{-1in}
                  \setlength{\voffset}{-1in}
                  \addtolength{\hoffset}{#5}
                            \addtolength{\voffset}{#6}
                            \setlength{\oddsidemargin}{#1}
                            \setlength{\evensidemargin}{#2}
                            \setlength{\textwidth}{\paperwidth}
                  \addtolength{\textwidth}{-\oddsidemargin}
                  \addtolength{\textwidth}{-\evensidemargin}
                  \addtolength{\textwidth}{-\marginparsep}
                  \addtolength{\textwidth}{-\marginparwidth}
                            \setlength{\topmargin}{#3}
                            \setlength{\textheight}{\paperheight}
                  \addtolength{\textheight}{-\topmargin}
                  \addtolength{\textheight}{-\headheight}
                  \addtolength{\textheight}{-\headsep}
                  \addtolength{\textheight}{-\footskip}
                  \addtolength{\textheight}{-#4}}
\pageformat{2cm}{3cm}{25mm}{25mm}{1pt}{0pt}

\usepackage{ifthen}
\newboolean{article}
    \setboolean{article}{true}
\newboolean{report}
\newboolean{book}
\newboolean{letter}
\newboolean{german}
\newboolean{italian}
\newboolean{nobaselinestretch}
\newboolean{nosectionappendix}
\newboolean{oldtoc}
\newboolean{nosectionequn}
\newboolean{notheorem}

\ifthenelse{\boolean{german}}{
    \usepackage{german}}{}

\usepackage[latin1]{inputenc}

\usepackage{amsmath}
\usepackage{amssymb}
\usepackage[mathscr]{eucal}

\ifthenelse{\boolean{notheorem}}{}{
    \usepackage{theorem}}

\ifthenelse{\boolean{nobaselinestretch}}{}{
    \renewcommand{\baselinestretch}{1.25}}

\newenvironment{env}[2]{\begin{#1}#2\end{#1}}{}
    \newcommand{\beq}[1]{\begin{env}{equation}{#1}}
    \newcommand{\beqn}[1]{\begin{env}{equation*}{#1}}
    \newcommand{\bal}[1]{\begin{env}{align}{#1}}
    \newcommand{\baln}[1]{\begin{env}{align*}{#1}}
    \newcommand{\bga}[1]{\begin{env}{gather}{#1}}
    \newcommand{\bgan}[1]{\begin{env}{gather*}{#1}}
    \newcommand{\bflal}[1]{\begin{env}{flalign}{#1}}
    \newcommand{\bflaln}[1]{\begin{env}{flalign*}{#1}}
    \newcommand{\bmu}[1]{\begin{env}{multline}{#1}}
    \newcommand{\bmun}[1]{\begin{env}{multline*}{#1}}
    \newcommand{\bsp}[1]{\begin{env}{split}{#1}}

    \newcommand{\eeq}{\end{env}}
    \newcommand{\eeqn}{\end{env}}
    \newcommand{\eal}{\end{env}}
    \newcommand{\ealn}{\end{env}}
    \newcommand{\ega}{\end{env}}
    \newcommand{\egan}{\end{env}}
    \newcommand{\eflal}{\end{env}}
    \newcommand{\eflaln}{\end{env}}
    \newcommand{\emu}{\end{env}}
    \newcommand{\emun}{\end{env}}
    \newcommand{\esp}{\end{env}}

\newcommand{\lf}{\vspace{2ex}}

\renewcommand{\bf}[1]{\textbf{#1}}
\renewcommand{\it}[1]{\textit{#1}}

\renewcommand{\sf}[1]{\textsf{#1}}

\renewcommand{\tt}[1]{\texttt{#1}}
\newcommand{\hl}[1]{\bf{\it{#1}}}

\newcommand{\mbf}[1]{\mathbf{#1}}
\newcommand{\msf}[1]{\text{\small$\sf{#1}$}}

\newcommand{\cmc}[1]{\mathcal{#1}}
\newcommand{\eus}[1]{\mathscr{#1}}
\newcommand{\euf}[1]{\mathfrak{#1}}
\newcommand{\bb}[1]{\mathbb{#1}}

\newcommand{\nbd}[1]{$#1$\nobreakdash--}
\newcommand{\ol}[1]{\overline{#1}}

\newcommand{\wh}[1]{\widehat{#1}}
\newcommand{\ve}{\varepsilon}
\newcommand{\vt}{\vartheta}

\newcommand{\abs}[1]{\left\lvert#1\right\rvert}
\newcommand{\norm}[1]{\left\lVert#1\right\rVert}
\newcommand{\babs}[1]{\bigl\lvert#1\bigr\rvert}

\newcommand{\snorm}[1]{\norm{\smash{#1}}}

\newcommand{\family}[1]{\left(#1\right)}

\newcommand{\bfam}[1]{\bigl(#1\bigr)}

\newcommand{\AB}[1]{\langle#1\rangle}

\newcommand{\CB}[1]{\{#1\}}
\newcommand{\bCB}[1]{\bigl\{#1\bigr\}}
\newcommand{\BCB}[1]{\Bigl\{#1\Bigr\}}
\newcommand{\SB}[1]{[#1]}
\newcommand{\bSB}[1]{\bigl[#1\bigr]}

\newcommand{\RO}[1]{[#1)}

\newcommand{\set}[2][]{
    \ifthenelse{\equal{#1}{}}{
        \CB{#2}}{
        \CB{#1~|~#2}}}
\newcommand{\bset}[2][]{
    \ifthenelse{\equal{#1}{}}{
        \bCB{#2}}{
        \bCB{#1~|~#2}}}
\newcommand{\Bset}[2][]{
    \ifthenelse{\equal{#1}{}}{
        \BCB{#2}}{
        \BCB{#1~\big|~#2}}}

\DeclareMathOperator{\ls}{\normalfont\msf{span}}
\DeclareMathOperator{\cls}{\ol{\ls}}

\DeclareMathOperator{\id}{\normalfont\msf{id}}

\newcommand{\N}{\bb{N}}

\newcommand{\R}{\bb{R}}
\newcommand{\bS}{\bb{S}}

\newcommand{\cB}{\cmc{B}}

\newcommand{\sB}{\eus{B}}

\newcommand{\sK}{\eus{K}}

\newcommand{\sN}{\eus{N}}

\newcommand{\sS}{\eus{S}}

\newcommand{\eR}{\euf{R}}

\newcommand{\U}{\mbf{1}}

\newcommand{\I}{{I\!\!\!\;I}}

\ifthenelse{\boolean{nosectionequn}}{}{
    \numberwithin{equation}{section}
    }

\ifthenelse{\boolean{article}\or\boolean{letter}\or\boolean{nosectionequn}}{
    \setboolean{nosectionappendix}{true}}{}
\ifthenelse{\boolean{nosectionappendix}}{}{
    \renewcommand{\appendix}{
        \chapter*{\appendixname}
        \addcontentsline{toc}{chapter}{\appendixname}
        \renewcommand{\thesection}{\Alph{section}}
        \setcounter{section}{0}}}
   
\ifthenelse{\boolean{report}\or\boolean{book}}{
    }{}

\ifthenelse{\boolean{notheorem}}{}{
        \newcommand{\mnname}{Mathematical note.}
        \newcommand{\enname}{End of the note.}
        \newcommand{\definame}{Definition.}
        \newcommand{\propname}{Proposition.}
        \newcommand{\lemname}{Lemma.}
        \newcommand{\exname}{Example.}
        \newcommand{\exername}{Exercise.}
        \newcommand{\remname}{Remark.}
        \newcommand{\obname}{Observation.}
        \newcommand{\thmname}{Theorem.}
        \newcommand{\corname}{Corollary.}
        \newcommand{\proofname}{Proof.}
    \ifthenelse{\boolean{german}}{
        \renewcommand{\mnname}{Mathematische Notiz.}
        \renewcommand{\enname}{Ende der Notiz.}
        \renewcommand{\exname}{Beispiel.}
        \renewcommand{\exername}{Übung.}
        \renewcommand{\remname}{Bemerkung.}
        \renewcommand{\obname}{Beobachtung.}
        \renewcommand{\thmname}{Satz.}
        \renewcommand{\corname}{Korollar.}
        \renewcommand{\proofname}{Beweis.}}{}
    \ifthenelse{\boolean{italian}}{
        \renewcommand{\mnname}{Nota matematica.}
        \renewcommand{\enname}{Fina della nota.}
        \renewcommand{\definame}{Definizione.}
        \renewcommand{\propname}{Proposizione.}
        \renewcommand{\exname}{Esempio.}
        \renewcommand{\exername}{Esercizio.}
        \renewcommand{\remname}{Nota.}
        \renewcommand{\obname}{Osservazione.}
        \renewcommand{\thmname}{Teorema.}
        \renewcommand{\corname}{Corollario.}
        \renewcommand{\proofname}{Dimostrazione.}

       \renewcommand{\appendixname}{Appendice}

       }{}
    \theoremheaderfont{\normalfont\bfseries}
    \theoremstyle{change}
        \theorembodyfont{\rmfamily}
            \newtheorem{emp}{}[section]
                \newcommand{\bemp}[1][]{
                    \begin{emp}\hskip-\labelsep\bf{#1}\hskip\labelsep}
                \newcommand{\eemp}{\end{emp}}
\newtheorem{itemp}[emp]{}
                \newcommand{\bitemp}[1][]{
                    \begin{itemp}\hskip-\labelsep\bf{#1}\hskip\labelsep\normalfont\itshape}
                \newcommand{\eitemp}{\end{itemp}}
            \newtheorem{mn}[emp]{\mnname}
                \newcommand{\bnm}{\begin{mn}~\begin{quotation}\renewcommand{\baselinestretch}{1}\small\noindent\ignorespaces}
                \newcommand{\enm}{\end{quotation}\hfill\bf{\enname}\end{mn}}
            \newtheorem{ex}[emp]{\exname}
                \newcommand{\bex}{\begin{ex}}
                \newcommand{\eex}{\end{ex}}
            \newtheorem{exer}[emp]{\exername}
                \newcommand{\bexer}{\begin{exer}}
                \newcommand{\eexer}{\end{exer}}
            \newtheorem{defi}[emp]{\definame}
                \newcommand{\bdefi}{\begin{defi}}
                \newcommand{\edefi}{\end{defi}}
            \newtheorem{rem}[emp]{\remname}
                \newcommand{\brem}{\begin{rem}}
                \newcommand{\erem}{\end{rem}}
            \newtheorem{ob}[emp]{\obname}
                \newcommand{\bob}{\begin{ob}}
                \newcommand{\eob}{\end{ob}}
        \theorembodyfont{\normalfont\itshape}
            \newtheorem{thm}[emp]{\thmname}
                \newcommand{\bthm}{\begin{thm}}
                \newcommand{\ethm}{\end{thm}}
            \newtheorem{prop}[emp]{\propname}
                \newcommand{\bprop}{\begin{prop}}
                \newcommand{\eprop}{\end{prop}}
            \newtheorem{cor}[emp]{\corname}
                \newcommand{\bcor}{\begin{cor}}
                \newcommand{\ecor}{\end{cor}}
            \newtheorem{lem}[emp]{\lemname}
                \newcommand{\blem}{\begin{lem}}
                \newcommand{\elem}{\end{lem}}
\newenvironment{empn}[1]{\lf\noindent\bf{#1}\ignorespaces\hskip\labelsep}{\lf}
		\newcommand{\bempn}[1]{\begin{empn}{#1}}
		\newcommand{\eempn}{\end{empn}}
		\newcommand{\bitempn}[1]{\begin{empn}{#1}\normalfont\itshape}
		\newcommand{\eitempn}{\end{empn}}
                \newcommand{\bnmn}{\begin{empn}{\mnname}~\begin{quotation}\renewcommand{\baselinestretch}{1}\small\noindent\ignorespaces}
                \newcommand{\enmn}{\end{quotation}\hfill\bf{\enname}\end{empn}}
		\newcommand{\bexn}{\begin{empn}{\exname}}
		\newcommand{\eexn}{\end{empn}}
		\newcommand{\bexern}{\begin{empn}{\exername}}
		\newcommand{\eexern}{\end{empn}}
		\newcommand{\bdefin}{\begin{empn}{\definame}}
		\newcommand{\edefin}{\end{empn}}
		\newcommand{\bremn}{\begin{empn}{\remname}}
		\newcommand{\eremn}{\end{empn}}
		\newcommand{\bobn}{\begin{empn}{\obname}}
		\newcommand{\eobn}{\end{empn}}

\newcommand{\qedsymbol}{~\rule[-0.35mm]{2mm}{2mm}}
    \newcounter{proof}[emp]
    \newenvironment{Proof}[1]{
        \vspace{1ex}
        \renewcommand{\item}[1][\stepcounter{proof}(\roman{proof})]%
            {##1\hskip\labelsep}
        \noindent\textsc{#1\hskip\labelsep}}{
        \nolinebreak\qedsymbol}
    \newcommand{\proof}[1][\proofname]{
        \begin{Proof}{#1}\ignorespaces}
    \newcommand{\qed}{\end{Proof}}
    \newcommand{\noqed}{
        \renewcommand{\qedsymbol}{}
        \end{Proof}}}
    \ifthenelse{\boolean{italian}}{
        \renewcommand{\proofname}{Dimostrazione.}}{}

\usepackage[varg]{txfonts}

\begin{document}

\title{$E_0$--Semigroups for Continuous Product Systems}
\author{}
\author{
~\\
Michael Skeide\thanks{This work is supported by research funds of University of Molise and Italian MIUR.}\\\\
{\small\itshape Universit\`a\ degli Studi del Molise}\\
{\small\itshape Dipartimento S.E.G.e S.}\\
{\small\itshape Via de Sanctis}\\
{\small\itshape 86100 Campobasso, Italy}\\
{\small{\itshape E-mail: \tt{skeide@unimol.it}}}\\
{\small{\itshape Homepage: \tt{http://www.math.tu-cottbus.de/INSTITUT/lswas/\_skeide.html}}}\\
}
\date{July 2006}

{
\renewcommand{\baselinestretch}{1}
\maketitle


\begin{abstract}
\noindent
We show that every \it{continuous} product system of correspondences over a \it{unital} \nbd{C^*}algebra occurs as the product system of a \it{strictly continuous} \nbd{E_0}semigroup.
\end{abstract}

}


{\parskip0.5ex plus 0.5ex minus 0.5ex

\section{Introduction}

\nbd{E_0}Semigroups on $\sB^a(E)$, the algebra of all adjointable operators on a Hilbert module $E$ over a \nbd{C^*}algebra $\cB$, give rise to product systems of correspondences over $\cB$. The first construction of this sort is due to Arveson in his trailblazing paper \cite{Arv89} which marks the begin of the modern theory of product systems. It took a whole serious of papers (Arveson \cite{Arv89,Arv90a,Arv89a,Arv90}) to show the converse statement, namely, that every product system of Hilbert spaces (\hl{Arveson system}, in the sequel) arises as the product system of an \nbd{E_0}semigroup on $\sB(H)$ for a Hilbert space $H$. For a long time there were no other proofs of this fact. Recently Liebscher \cite{Lie00p1} provided a different still very involved proof. In Skeide \cite{Ske06} we found a short and self-contained proof and shortly after Arveson \cite{Arv06} presented yet another short proof. In Skeide \cite{Ske06p4} we showed that Arveson's construction in \cite{Arv06} leads to a result that is unitarily equivalent to a special case of the construction in \cite{Ske06}.

The first construction of a product system from an \nbd{E_0}semigroup on a general $\sB^a(E)$ is done in Skeide \cite{Ske02} under the assumption that $E$ has a \it{unit vector}.\footnote{Apparently, there is a construction of a product systems from \nbd{E_0}semigroups on type II$_1$ factors due to Alevras in his thesis. But, still after several inquiries we do not have this thesis available.} The general case for nonunital \nbd{C^*}algebras (that is, in particular, without unit vectors) is discussed in Skeide \cite{Ske04p} (based on the representation theory of $\sB^a(E)$ in Muhly, Skeide and Solel \cite{MSS06}).

It is the scope of these notes to prove the converse, every product system comes from an \nbd{E_0}semigroup, in the special case of \it{continuous} product systems of correspondences over a unital \nbd{C^*}algebra. The general case will be treated in Skeide \cite{Ske06p1}. In several places in these notes (Remark \ref{irem}, Section \ref{concSEC}) we will comment on what the differences and additional complications are (beginning with a whole bunch of technically quite different variants). What made us decide to publish the present case separately are two reasons: Firstly, it allows a complete solution of the problem. By this we mean that we have a complete correspondence between a sufficiently interesting class of \nbd{E_0}semigroups on the one side and a handy class of product systems on the other side. Secondly, the special properties allow for a particularly simple treatment, immitating the construction in \cite{Arv06}. (Anyway, we point out in Section \ref{concSEC} that our conditions are not that special. In fact, that part of the condition that allows to apply Arveson's construction are fulfilled by every Arveson system.)

\section{The product system of an $E_0$--semigroup}\label{apsSEC}

By $\bS$ we denote either $\N_0=\CB{0,1,\ldots}$ (\hl{discrete} case) or $\R_+=\RO{0,\infty}$ (\hl{continuous time} case). Let $E$ be a Hilbert \nbd{\cB}module. Suppose that $\vt=\bfam{\vt_t}_{t\in\bS}$ is an \hl{\nbd{E_0}semigroup} on $\sB^a(E)$. By this we mean that $\vt$ is a semigroup of unital endomorphisms $\vt_t$ of $\sB^a(E)$. In these notes we shall always assume that the $\vt_t$ are \hl{strict} (that is, they are continuous for the strict topology on bounded subsets of $\sB^a(E)$ or, equivalently, already the action of $\vt_t(\sK(E))$ of the compact operators $\sK(E)$ on $E$ via $\vt_t$ is nondegenerate).
}

As discussed in \cite{Ske04p}, using the results from \cite{MSS06}, with every $\vt_t$ $(t>0)$ we may associate a correspondence $E_t$ over $\cB$ and a unitary $v_t\colon E\odot E_t\rightarrow E$ such that $\vt_t(a)=\vt^v_t(a):=v_t(a\odot\id_t)v_t^*$. Moreover, there are bilinear unitaries $u_{s,t}\colon E_s\odot E_t\rightarrow E_{s+t}$ such that $(x_ry_s)z_t=x_r(y_sz_t)$ where, following Arveson's convention, we denote $x_sy_t:=u_{s,t}(x_s\odot y_t)$. Putting $E_0=\cB$ (the \it{trivial} correspondence over $\cB$) the families $v_t$ and $u_{s,t}$ extend to time $0$ by the canonical idenitifications. In other words, the family $E^\odot=\bfam{E_t}_{t\in\bS}$ is a \hl{product system} in the sense of Bhat and Skeide \cite[Definition 4.7]{BhSk00}. If we wish to underline absence of continuity or measurability conditions, we say $E^\odot$ is an \hl{algebraic} product system. Also, using the same notation $xy_t:=v_t(x\odot y_t)$, the $v_t$ fulfill $(xy_s)z_t=x(y_sz_t)$.

We do not give details as we shall discuss below a different construction in the special case when $E$ has a \hl{unit vector} $\xi$ (that is, $\AB{\xi\xi}=\U$ so that, in particular, $\cB$ is unital). We shall say, $E$ is \hl{unital}, if it admits a unit vector. We just mention that, actually, all $E_t$ may be viewed as correspondences over the \hl{range ideal} $\cB_E:=\cls\AB{E,E}$ of $E$. (For the right action this is trivial as every Hilbert \nbd{\cB}module may be considered as a Hilbert \nbd{\cB_E}module. For the left action it means that already the action of $\cB_E$ on $E$ is nondegenerate.) Furthermore, $E^\odot$ together with a family $v_t$ such that $\vt^v_t$ gives back $\vt_t$ is determined up to isomorphism of product systems by these properties. See \cite{Ske04p} for details. So, replacing $\cB$ (and, therefore, also $E_0$) with $\cB_E$ we may assume that $E$ is \hl{full}, that is, $\cB=\cB_E$. In this case, the constructions give the correct structures also for time $0$. Also, all $E_t$ are necessarily full. In general, for full product system $E^\odot$ (that is, all $E_t$ are full) by a \hl{left dilation} of $E^\odot$ to a full Hilbert module $E$ we shall understand a family of unitaries $v_t\colon E\odot E_t\rightarrow E$ that fulfill $(xy_s)z_t=x(y_sz_t)$. (Note that, using associativity, $v_0$ is bound to be the canonical identification. In \cite{Ske06p1} we will discuss a version that includes also the nonfull case.) For every left dilation the $\vt^v_t$ define an \nbd{E_0}semigroup on $\sB^a(E)$. So, finding an \nbd{E_0}semigroup for a full product system is equivalent to find a left dilation.

We will now discuss the construction from Skeide \cite{Ske02} of an algebraic product system $E^\odot$ from an \nbd{E_0}semigroup $\vt$ on $\sB^a(E)$ for unital $E$. (This is a direct generalization of Bhat's construction in \cite{Bha96} of an Arveson system from an \nbd{E_0}semigroup on $\sB(H)$.) While in Section \ref{cpsSEC} we will assume that $\vt$ is \hl{strictly continuous} (that is, for every $a\in\sB^a(E)$ the function $t\mapsto\vt_t(a)$ is strictly continuos) or, equivalently (because all $\vt_t$ are contractive \nbd{*}maps), strongly continuous. It is the unit vector which will allow us, as in Skeide \cite{Ske03b}, to define a continuous structure on $E^\odot$.

\brem\label{unifullrem}
As long as $E$ is full over a unital \nbd{C^*}algebra (or slightly more weakly, as long as $\cB_E$ is unital) the assumption of a unit vector is not critical. Indeed, \cite[Lemma 2.2]{Ske04p} asserts that a finite direct sum $E^n$ of copies of $E$ has a unit vector. By inflation the \nbd{E_0}semigroup $\vt$ on gives rise to an \nbd{E_0}semigroup $\vt^n$ on $\sB^a(E^n)=M_n(\sB^a(E))$ and $\vt^n$ is strictly continuous, if and only if $\vt$ is strictly continuous. It is not difficult to show that $\vt$ and $\vt^n$ have the same product system. (This follows simply because a left dilation $v_t$ of $E^\odot$ to $E$ gives rise to a left dilation $v^n_t\colon E^n\odot E_t\rightarrow E^n$ which induces $\vt^n$ and from the uniqueness of the product system.) In the case when $\cB_E$ is nonunital, so that it is meaningless to ask for a unit vector, we do not know how to impose a continuos structure on the product system $E^\odot$. Anyway, the left dilation we are going to construct from a continuous product system will be to a unital Hilbert module, so in our context it is perfectly admissible to restrict our considerations to left dilations to unital Hilbert modules.
\erem

If $\xi\in E$ is a unit vector, we put $E_t:=\vt_t(\xi\xi^*)E$. On $E_t$ we define a left action of $\cB$ by setting $bx_t:=\vt_t(\xi b\xi^*)x_t$. This left action is unital, so that $E_t$ is a correspondence over $\cB$. It is easy to check that
\beqn{
x\odot y_t
~\longmapsto~
\vt_t(x\xi^*)y_t
}\eeqn
defines an isometry $v_t\colon E\odot E_t\rightarrow E$. Surjectivity follows from strictness of $\vt_t$; see \cite{Ske02,Ske03b} for details. Obviously, $\vt_t(a)=v_t(a\odot\id_t)v_t^*$. The restriction of $v_t$ to $E_s\odot E_t\subset E\odot E_t$ defines a bilinear unitary $u_{s,t}$ onto $E_{s+t}$. (Clearly, $u_{s,t}$ is into $E_{s+t}$ and bilinear. Surjectivity follows from $\vt_{s+t}(\xi\xi^*)v_t=v_t(\vt_s(\xi\xi^*)\odot\id_t)$.) We leave it as an exercise to check that $E^\odot:=\bfam{E_t}_{t\in\bS}$ with the $u_{s,t}$ is a (full, of course) product system and the the $v_t$ form a left dilation of $E^\odot$ to $E$; see again \cite{Ske02,Ske03b} for details.

By \cite[Proposition 2.3]{Ske02} the product system does not depend on the choice of the unit vector. Indeed, if $\xi'$ is another unit vector, then $\vt_t(\xi'\xi^*)$ defines an isomorphism from the product system constructed from $\xi$ to the product system constructed from $\xi'$.

We close this section on algebraic product systems with the construction from \cite{BhSk00} (\cite{BBLS04} for the general case) of an \nbd{E_0}semi\-group for a product system when this product system has a unital unit. A \hl{unit} in a product system $E^\odot=\bfam{E_t}_{t\in\bS}$ of correspondences $E_t$ over a unital \nbd{C^*}algebra $\cB$ is a family $\xi^\odot=\bfam{\xi_t}_{t\in\bS}$ of vectors $\xi_t\in E_t$ with $\xi_0=\U$ which compose as $\xi_s\xi_t=\xi_{s+t}$. The unit is \hl{unital}, if every $\xi_t$ is a unit vector. (If a unit $\xi^\odot$ is \hl{continuous} in the sense that the CP-semigroup $T=\bfam{T_t}_{t\in\bS}$ on $\cB$ defined by setting $T_t(b)=\AB{\xi_t,b\xi_t}$ is uniformly continuous, then this unit may be ``normalized'' to a unital unit within the product system $E^\odot$; see Skeide and Liebscher \cite[Example 4.2]{LiSk05p}.)

We observe that a unital unit gives rise to an inductive system of isometric embeddings $\xi_s\odot\id_t\colon E_t\rightarrow E_s\odot E_t$ defined by setting $(\xi_s\odot\id_t)x_t=\xi_sx_t$.

\bitemp[Theorem {\protect\cite[Section 4.4]{BBLS04}}.~]\label{uniE0thm}
Let $E$ denote the inductive limit over $E_t$. All $\xi_t\in E_t$ are imbedded to the same unit vector in $E$ which we denote by $\xi$. For every $t\in\bS$ the factorization $u_{s,t}\colon E_s\odot E_t\rightarrow E_{s+t}$ for $s\to\infty$ gives rise to a factorization $v_t\colon E\odot E_t\rightarrow E$ of the inductive limit and the $v_t$ form a left dilation of $E^\odot$ to $E$. Moreover, $\xi\xi_t=\xi$ so that the product system of the \nbd{E_0}semigroup $\vt$ on $\sB^a(E)$ defined by setting $\vt_t(a)=v_t(a\odot\id_t)v_t^*$ is $\vt_t(\xi\xi^*)E=E_t$ including the correct product system structure.
\eitemp

\section{Continuous product systems}\label{cpsSEC}

We pass now to the continuous case. The following Definition of continuous product system is \cite[Definition 7.1]{Ske03b} except that we have removed that $\cB$ is assumed unital. It is motivated by the fact that every strictly continuous \nbd{E_0}semigroup acting on the operators of a unital Hilbert module fulfills these requirements.

\bdefi\label{cPSdef}
Let $E^\odot=\bfam{E_t}_{t\in\R_+}$ be a product system of correspondences over a \nbd{C^*}algebra $\cB$ with a family $i=\bfam{i_t}_{t\in\R_+}$ of isometric embeddings $i_t\colon E_t\rightarrow\wh{E}$ into a Hilbert \nbd{\cB}module $\wh{E}$. Denote by
$$
CS_i(E^\odot)
~=~
\BCB{x=\bfam{x_t}_{t\in\R_+}\colon x_t\in E_t,t\mapsto i_tx_t\mbox{~is continuous}}
$$
the set of \hl{continuous sections} of $E^\odot$ (with respect to $i$). We say $E^\odot$ is \hl{continuous} (with respect to $i$), if the following conditions are satisfied.
\begin{enumerate}
\item\label{csdef}
For every $y_t\in E_t$ we can find a continuous section $x\in CS_i(E^\odot)$ such that $x_t=y_t$.

\item \label{cpsdef}
For every pair $x,y\in CS_i(E^\odot)$ of continuous sections the function
$$
(s,t)
~\longmapsto~
i_{s+t}(x_sy_t)
$$
is continuous.
\end{enumerate}
We say two embeddings $i$ and $i'$ have the same \hl{continuous structure}, if $CS_i(E^\odot)=CS_{i'}(E^\odot)$.
\edefi

Roughly speaking, $E^\odot$ is a Banach subbundle of the trivial Banach bundle $\R_+\times\wh{E}$ that contains enough continuous sections and the product system structure respects continuity of sections. Note also that by \cite[Proposition 7.9]{Ske03b} Condition 1 may be replaced by the weaker condition that for every $t\in\R_+$ the set $\CB{x_t\colon x\in CS_i(E^\odot)}$ is total in $E_t$. (The proof does not depend on that the definition here is slightly more general, and presents a typical example of dealing with continuous sections.) Note also that Condition 2 may be replaced with the weaker condition that the function $(s,t)\mapsto \AB{z,i_{s+t}(x_sy_t)}$ is continuous for every $z\in\wh{E}$ and every pair $x,y\in CS_i(E^\odot)$. (The proof is very much analogue to that of the well-known fact that on the unitary group of a Hilbert space strong and weak topology coincide.)

Observe that, in particular, for every $z_s\in E_s$ and every section $x\in CS_i(E^\odot)$ the functions $t\mapsto i_{s+t}(z_sx_t)$ and $t\mapsto i_{t+s}(x_tz_s)$ are continuous. (Simply choose a section $y\in CS_i(E^\odot)$ with $y_s=z_s$ and keep $s$ constant.)

Before we investigate the continuous structure of the product system of a strictly continuous \nbd{E_0}semigroup we illustrate how strong the condition to be continuous at $t=0$ is for a product system.

\blem\label{cflem}
If $\cB$ is unital, then a continuous product system $E^\odot$ of correspondences over $\cB$ contains a continuous section $\zeta\in CS_i(E^\odot)$ that consists entirely of unit vectors and fulfills $\zeta_0=\U$. In particular, every $E_t$ contains a unit vector (and, therefore, is full).
\elem

\proof
By assumption \ref{cPSdef}\eqref{csdef} there exists a continuous section $x$ such that $x_0=\U\in\cB=E_0$. As $\U$ is invertible and the invertible elements form an open subset of $\cB$, the elements $\abs{x_t}:=\sqrt{\AB{x_t,x_t}}$ are invertible on an interval $\SB{0,\ve}$ for a suitable $\ve>0$. We put $y_t=x_t\abs{x_t}^{-1}$ for $t\in\SB{0,\ve}$ and $y_t=x_t\abs{x_\ve}^{-1}$ for $t>\ve$. Clearly, this defines a continuous section $y$. We define a section $\zeta=\bfam{\zeta_t}_{t\in\R_+}$ by setting $\zeta_t=y_{t-n\ve}y_\ve^n$ where $n=\bSB{\frac{t}{\ve}}$ is the unique integer such that $t-n\ve\in\RO{0,\ve}$. By construction all $\zeta_t$ are unit vectors. As $y$ is continuous, the section $\zeta$ is continuous for all $t\notin\N_0\ve$. If $t=n\ve$, then the left and right limit limit at $t$ are
\baln{
\lim_{\delta\to+0}i_{t+\delta}\zeta_{t+\delta}
&
~=~
i_{t}(y_0y_\ve^n)
~=~
i_{t}y_\ve^n
~=~
i_t\zeta_t
\\
\lim_{\delta\to+0}i_{t-\delta}\zeta_{t-\delta}
&
~=~
i_{t}(y_\ve y_\ve^{n-1})
~=~
i_{t}y_\ve^n
~=~
i_t\zeta_t.
}\ealn
So $\zeta$ is continuous.\qed

\lf
Now let $E$ be a unital Hilbert \nbd{\cB}module and fix a unit vector $\xi\in E$. Suppose that $\vt=\bfam{\vt_t}_{t\in\R_+}$ is a strictly continuous \nbd{E_0}semigroup and dennote by $E^\odot$ its associated product system constructed rom $\xi$ as $E_t:=\vt_t(\xi\xi^*)E\subset E$. Put $\wh{E}:=E$ and let $i_t$ denote the canonical embeddings $E_t\rightarrow\wh{E}$. Choose $y_t\in E_t$. Then $x$ with $x_s:=\vt_s(\xi\xi^*)y_t$ is a continuous section such that $x_t=y_t$. Moreover, if $x,y\in CS_i(E^\odot)$ is a pair of continuous sections, then
\beqn{
(s,t)
~\longmapsto~
i_{s+t}(x_sy_t)
~=~
\vt_t(x_s\xi^*)y_t
}\eeqn
is, clearly, continuous.

So far, this has been explained in \cite[Section 7]{Ske03b}. But, we mention that if $\xi'\in E$ is another unit vector, then the isomorphism $\vt_t(\xi'\xi^*)$, clearly, sends continuous sections to continuous sections. Therefore, the continuous structures of the product system constructed from $\xi$ and of the product system constructed from $\xi'$ coincide.

Without proof we state the following which improves on \cite[Theorem 7.5]{Ske03b} where the unit $\xi^\odot$ was required to be continuous.

\bthm\label{cucthm}
Let $E^\odot$ be a continuous product system and let $\xi^\odot$ be a unital unit in $CS_i(E^\odot)$. Then the \nbd{E_0}semigroup $\vt$ constructed from $\xi^\odot$ by Theorem \ref{uniE0thm} is strictly continuous and the continuous structure derived from $\vt$ coincides gives back the continuous structure of $E^\odot$.
\ethm

The proof is the same as in \cite{Ske03b} but there we were interested only in continuous units. The scope was rather to start with an algebraic product system and a continuous unit, leading to an \nbd{E_0}semigroup by Theorem \ref{uniE0thm} which shows to be strictly continuous. We wanted to convince ourselves that the induced continuous structure does not depend on the unit as long as the units are sufficiently contionuous \it{with respect to eachother}.

Of course, also Theorem \ref{cucthm} shows that the continuous structure induced by a unital unit does not depend on the choice. It is interesting to note that, so far, we do not know whether the inductive limits contructed from different units are isomorphic. In fact, we strongly suspect that they need not be isomorphic.

Finally, we mention that the continuous structure of a product system associated with a strictly continuous \nbd{E_0}semigroup may equally well be expressed in terms of the left dilation that gives back the \nbd{E_0}semigroup. In fact, the canonical embedding $E_t=\vt_t(\xi\xi^*)\rightarrow E$ is nothing but $v_t(\xi\odot\id_t)\colon x_t\mapsto v_t(\xi\odot x_t)=\xi x_t$.

\section{The construction}\label{constrSEC}

\brem\label{irem}
The basic idea of Skeide \cite{Ske06} (which we describe immediately for modules) to find a left dilation of a product system $E^\odot=\bfam{E_t}_{t\in\R_+}$ was to start with a left dilation of the discrete subsystem $\bfam{E_t}_{t\in\N_0}$ to a Hilbert module $\breve{E}$, that is, with a family of unitaries $\breve{v}_n\colon\breve{E}\odot E_n\rightarrow\breve{E}$ that fulfill the necessary associativity conditions. We put $E:=\breve{E}\odot\int_0^1 E_\alpha\,d\alpha$. The following identifications
\bal{
\notag
E\odot E_t
~=~
\breve{E}\odot\family{\int_0^1E_\alpha\,d\alpha}\odot E_t
&
~=~
\breve{E}\odot\int_t^{1+t}E_\alpha\,d\alpha
\\[2ex]
\notag
&
~\cong~
\family{\breve{E}\odot E_n\odot\int_{t-n}^1E_\alpha\,d\alpha}
\oplus
\family{\breve{E}\odot E_{n+1}\odot\int_0^{t-n}E_\alpha\,d\alpha}
\\[2ex]\label{idea}
&
~\cong~
\family{\breve{E}\odot\int_{t-n}^1E_\alpha\,d\alpha}
\oplus
\family{\breve{E}\odot\int_0^{t-n}E_\alpha\,d\alpha}
~=~
E
}\eal
suggest, then, a family of unitaries $v_t\colon E\odot E_t\rightarrow E$. The slightly tedious thing in \cite{Ske06} was to show associativity, that, is that the $v_t$ form a dilation to $E$. But, by that method whenever we are able to dilate the discrete subsystem $\bfam{E_t}_{t\in\N_0}$ of $E^\odot$ we are also able to dilate the whole product system $E^\odot$.

Existence of the dilation of the discrete subsystem was settled in \cite[Theorem 6.6]{Ske04p} for full correspondences over a unital \nbd{C^*}algebra and in \cite[main theorem]{Ske04p} for strongly full von Neumann correspondences. Here, for continuous product systems of full correspondences over a unital \nbd{C^*}algebra, the situation is even better. By Lemma \ref{cflem} $E_1$ contains a unit vector $\zeta_1$. We do not know whether $E^\odot$ has a unital unit. (In this case, the whole construction in the remainder would be superfluous, as we could simply apply Theorem \ref{uniE0thm}.) But, at least the discrete subsystem $\bfam{E_n}_{n\in\N_0}$ has a unital unit, namely, $\xi^\odot=\bfam{\xi_n}_{n\in\N_0}$ with $\xi_n:=\zeta_1^n$. So, Theorem \ref{uniE0thm} provides us with a dilation at least of the discrete subsystem we can use as input for the construction as indicated in \eqref{idea}.

It is this case, fixing a unit vector $\zeta_1\in E_1$, which was by treated by Arveson \cite{Arv06} in a different way. Roughly speaking, as explained in \cite{Ske06p4}, the construction of \cite{Arv06} can be interpreted as exchanging in $E=\breve{E}\odot\int_0^1 E_\alpha\,d\alpha$ the construction of the inductive limit (giving $\breve{E}$) and the direct integral, and giving a very concrete interpretation of the elements of the inductive limit over $E_n\odot\int_0^1 E_\alpha\,d\alpha=\int_n^{n+1}E_\alpha\,d\alpha$ in terms of sections of the product system with a handy equivalence relation. As we do have unit vectors, we follow the same construction here.
\erem

Now we start with the construction. But, before we can really start we have to say a few words about the direct integrals. If $E^\odot$ is a continuous product system with continuous structure defined by the family $i$ of embeddings $i_t\colon E_t\rightarrow\wh{E}$, then every section $x=\bfam{x_t}_{t\in\R_+}$ in $E^\odot$ gives rise to a function $t\mapsto x(t):=i_tx_t$ with values in $\wh{E}$. Let $0\le a<b<\infty$. By $\int_a^b E_\alpha\,d\alpha$ we understand the norm closure of the pre-Hilbert module that consists of continuous sections $x\in CS_i(E^\odot)$ restricted to $\RO{a,b}$ with inner product
\beqn{
\AB{x,y}_{\SB{a,b}}
~:=~
\int_a^b\AB{x_\alpha,y_\alpha}\,d\alpha
~=~
\int_a^b\AB{x(\alpha),y(\alpha)}\,d\alpha.
}\eeqn
Note that all continuous sections are bounded on the compact interval $\SB{a,b}$ and, therefore, square integrable. As by Lemma \ref{cflem} there is a continuous section of unit vectors, also $\int_a^b E_\alpha\,d\alpha$ contains a unit vector.

\bprop\label{rcllprop}
$\int_a^b E_\alpha\,d\alpha$ contains the space $\eR_{\RO{a,b}}$ of restrictions to $\RO{a,b}$ of those sections $x$ for which $t\mapsto x(t)$ is right continuous with finite jumps (this implies that there exists a left limit) in finitely many points and bounded on $\RO{a,b}$, as a pre-Hilbert submodule.
\eprop

\proof
It is sufficient to observe that we may construct a right continuous section with left limit in each point which has a determined jump in one specific point, $r$ say, and is continuous otherwise on $\RO{a,b}\backslash\CB{r}$. Simply choose a continuous section that assumes at $r$ the jumpsize and multiply it by a sequence of continuous functions that \nbd{L^2}approximates the indicator function of $\RO{r,b}$. Then the limit of this sequence has the desired property. Adding up a finite number of such functions we produce a function that has exactly the same jumps as $x$, so that the difference is continuous and, therefore, in $\eR_{\RO{a,b}}$.

Of course, the inner product is definite on right continuous functions, that is, $\eR_{\RO{a,b}}$ is indeed a pre-Hilbert module. (This would fail, if we considered the interval $\SB{a,b}$.)\qed

\lf
Let $\sS$ denote the right \nbd{\cB}module of all sections $x=\bfam{x_t}_{t\in\R_+}$ of $E^\odot$ which are \hl{locally $\eR$}, that is, for every $0\le a<b<\infty$ the restriction of $x$ to $\RO{a,b}$ is in $\eR_{\RO{a,b}}$, and which are \hl{stable} with respect to the unit vector $\zeta_1\in E_1$, that is, there exists an $\alpha_0\ge0$ such that
\beqn{
x_{\alpha+1}
~=~
\zeta_1x_{\alpha}
}\eeqn
for all $\alpha\ge\alpha_0$. By $\sN$ we denote the subspace of all sections in $\sS$ which are eventually $0$, that is, of all sections $x\in\sS$ for which there exists an $\alpha_0\ge0$ such that $x_\alpha=0$ for all $\alpha\ge\alpha_0$. A straightforward verification shows that
\beqn{
\AB{x,y}
~:=~
\lim_{m\to\infty}\int_m^{m+1}\AB{x(\alpha),y(\alpha)}\,d\alpha
}\eeqn
defines a semiinner product on $\sS$ and that $\AB{x,x}=0$ if and only if $x\in\sN$. Actually, we have
\beqn{
\AB{x,y}
~=~
\int_T^{T+1}\AB{x(\alpha),y(\alpha)}\,d\alpha
}\eeqn
for all sufficiently large $T>0$; see \cite[Lemma 2.1]{Arv06}. So, $\sS/\sN$ becomes a pre-Hilbert module with inner product $\AB{x+\sN,y+\sN}:=\AB{x,y}$. By $E$ we denote its completion. (By arguments similar to those in \cite{Ske06p4}, the $E$ here, indeed, is canonically isomorphic to the $E$ discussed in Remark \ref{irem}.)

\bprop\label{denseprop}
For every section $x$ and every $\alpha_0\ge0$ define the section $x^{\alpha_0}$ as
\beqn{
x^{\alpha_0}_\alpha
~:=~
\begin{cases}
0&\alpha<\alpha_0
\\
\zeta_1^nx_{\alpha-n}&\alpha\in\RO{\alpha_0+n,\alpha_0+n+1},n\in\N_0.
\end{cases}
}\eeqn
If $x$ is in $CS_i(E^\odot)$, then $x^{\alpha_0}$ is in $\sS$. Moreover, the set $\bCB{x^{\alpha_0}+\sN\colon x\in CS_i(E^\odot),\alpha_0\ge0}$ is a dense submodule of $E$.
\eprop

\proof
Of course, $x^{\alpha_0}$ is in $\sS$ whenever $x$ is continuous. So, let $y$ be a section in $\sS$ and choose $\alpha_0$ such that $y_{\alpha+1}=\zeta_1y_{\alpha}$ for every $\alpha\ge\alpha_0$. Then, $y=y^{\alpha_0}\mod\sN$. Now, if $x_n$ is a sequence of continuous sections that approximates $y$ in $\eR_{\RO{\alpha_0,\alpha_0+1}}$, then $x_n^{\alpha_0}$ approximates $y^{\alpha_0}$ in $\sS/\sN$. That is, the set is dense in $E$. Of course, the set of sections of the form $x^{\alpha_0}$ is invariant under right multiplication and modulo $\sN$ also under addition.\qed

\lf
Note that for the continuous section $\zeta$ of unit vectors from Lemma \ref{cflem} also the section $\zeta^0$ is continuous. (This follows just as in the proof of Lemma \ref{cflem} now with $\ve=1$.)

\bcor\label{univeccor}
$\xi:=\zeta^0+\sN$ is a unit vector in $E$.
\ecor

\brem
Observe also that $\zeta^0_n=\xi_n=(\zeta^0_1)^n$. So $\zeta^0$ interpolates the unit $\xi^\odot$ of the discrete subsystem in Remark \ref{irem}.
\erem

After these preparations it is completely plain to see that for every $t\in\R_+$ the map $x\odot y_t\mapsto xy_t$, where
\beqn{
(xy_t)_\alpha
~=~
\begin{cases}
x_{\alpha-t}y_t&\alpha\ge t,
\\
0&\text{else},
\end{cases}
}\eeqn
defines an isometry $v_t\colon E\odot E_t\rightarrow E$, and that these isometries iterate associatively.

So far we discussed that part of the construction that is immediate, once the idea of the module of stable sections and its inner product are understood. The remaining work, surjectivity of the $v_t$, continuity of the \nbd{E_0}semigroup and compatibility of the continuous structure arising from that \nbd{E_0}semigroup with the original one, require a certain ammount of technical work and cover the remainder of this section.

\bprop\label{surprop}
Each $v_t$ is surjective.
\eprop

\proof
By Proposition \ref{denseprop} it is sufficient to approximate every section of the form $x^{\alpha_0}$ with $x\in CS_i(E^\odot),\alpha_0\ge0$ in the (semi-)inner product of $\sS$ by finite sums of sections of the form $yz_t$ for $y\in\sS,z_t\in E_t$. As what the section does on the finite interval $\RO{0,t}$ is not important for the inner product, we may even assume that $\alpha_0\ge t$. And as in the proof of Proposition \ref{denseprop} the approximation can be done by approximating $z$ in $\eR_{\RO{\alpha_0,\alpha_0+1}}$ and then extending the restriction to $\RO{\alpha_0,\alpha_0+1}$ stably to the whole axis.

So let $\alpha_0\ge t$ and let $x$ be a continuous section. We will approximate the contininuous section $\alpha\mapsto x_\alpha$ uniformly on the compact interval $\SB{\alpha_0,\alpha_0+1}$ (and, therefore, in $L^2$) by finite sums over sections of the form $\alpha\mapsto y_{\alpha-t}z_t$. Choose $\ve>0$. For every $\beta\in\SB{\alpha_0,\alpha_0+1}$ choose $n^\beta\in\N,y_k^\beta\in E_{\beta-t},z_k^\beta\in E_t$ such that $\snorm{x_\beta-\sum_{k=1}^{n^\beta}y_k^\beta z_k^\beta}\le\frac{\ve}{2}$. Choose continuous sections $\bar{y}_k^\beta=\bfam{(\bar{y}_k^\beta)_\alpha}_{\alpha\in\R_+}\in CS_i(E^\odot)$ such that $(\bar{y}_k^\beta)_{\beta-t}=y_k^\beta$. For every $\beta$ chose the maximal interval $I_\beta\subset\SB{\alpha_0,\alpha_0+1}$ such that $\snorm{x_\alpha-\sum_{k=1}^{n^\beta}(\bar{y}_k^\beta)_{\alpha-t}z_k^\beta}<\ve$ for all $\alpha\in I_\beta$. Of course, $I_\beta$ contains $\beta$ and is open in $\SB{\alpha_0,\alpha_0+1}$, because it is the inverse image of an open set under a continuous function. In other words, the $I_\beta$ form an open cover of the compact set $\SB{\alpha_0,\alpha_0+1}$ so that we may choose a finite subcover determined, say, by $m$ values $\beta_1,\ldots,\beta_m\in\SB{\alpha_0,\alpha_0+1}$. By taking away from $I_{\beta_i}$ everything that is already contained in $I_{\beta_1}\cup\ldots\cup I_{\beta_{i-1}}$, we define a finite partition $I_i$ of $\SB{\alpha_0,\alpha_0+1}$. Taking away the point $\alpha_0+1$ and adjusting the endpoints of the $I_i$ suitably, we may assume that all $I_i$ are right open. Denote by $\I_i$ the indicator function of $I_i$. Then, restriction of the piecewise continuous section
\beqn{
\alpha
~\longmapsto~
\begin{cases}
0&\alpha<t
\\
\displaystyle\sum_{i=1}^m\sum_{k=1}^{n_{\beta_i}}(\bar{y}_k^{\beta_i})_{\alpha-t}z_k^{\beta_i}\I_i(\alpha))&\alpha\ge t
\end{cases}
}\eeqn
to $\RO{\alpha_0,\alpha_0+1}$ is in $\eR_{\RO{\alpha_0,\alpha_0+1}}$ and approximates $\alpha\mapsto z_\alpha$ uniformly on $\RO{\alpha_0,\alpha_0+1}$ up to $\ve$.\qed

\lf
So, the $v_t$ form a dilation of $E^\odot$ to $E$. To show that the associated \nbd{E_0}semigroup is continuous, we show first that the dilation is \hl{continuous} in the following sense.

\bprop\label{cdprop}
For every $x\in E$ and every continuous section $y\in CS_i(E^\odot)$ the function $t\mapsto xy_t$ is continuous.
\eprop

\proof
As $y$ is bounded locally uniformly, it is sufficient to show the statement for all $x$ from a dense subset of $E$. So suppose that $x$ (modulo $\sN$) is given by a section in $\sS$ of the form $z^{\alpha_0}$ for $z\in CS_i(E^\odot)$ and $\alpha_0\ge0$. To calculate $\norm{z^{\alpha_0}y_t-z^{\alpha_0}y_s}$ we have to integrate over $\alpha$ the values of $\abs{z^{\alpha_0}_{\alpha-t}y_t-z^{\alpha_0}_{\alpha-s}y_s}^2$ for $\alpha$ in any unit interval such that  $\alpha-t$ and $\alpha-s$ are not smaller than $\alpha_0$. So
\beq{\label{secnd}
\abs{z^{\alpha_0}y_t-z^{\alpha_0}y_s}^2
~=~
\int_d^{d+1}\abs{z^{\alpha_0}_{\alpha+s}y_t-z^{\alpha_0}_{\alpha+t}y_s}^2\,d\alpha
}\eeq
for all $d\ge\alpha_0$. The function $(\alpha,t)\mapsto z_\alpha y_t$ is uniformly continuous on each interval $\RO{\alpha_0+n,\alpha_0+n+1}\times\SB{a,b}$ and it is bounded on every $\R_+\times\SB{a,b}$. We fix a $t$, we choose a (sufficiently big) $d$ such that $n=d+t-\alpha_0$ is an integer and we choose $\ve\in\bfam{0,\frac{1}{2}}$. Then in
\beqn{
\abs{z^{\alpha_0}y_t-z^{\alpha_0}y_s}^2
\\
~=~
{\Biggl[\textstyle\int_d^{d+\ve}+\int_{d+\ve}^{d+1-\ve}+\int_{d+1-\ve}^{d+1}\Biggr]}~\babs{\smash{z^{\alpha_0}_{\alpha+s}y_t-z^{\alpha_0}_{\alpha+t}y_s}}^2\,d\alpha
}\eeqn
the first and the last integral are bounded by $\ve$ times a constant which is independent on $s\in\bfam{t-\frac{1}{2},t+\frac{1}{2}}$. For such $s$, in the middle integral both $\alpha-t$ and $\alpha-s$ are in the same interval $\RO{\alpha_0+n,\alpha_0+n+1}$, so that both $z^{\alpha_0}_{\alpha+s}y_t$ and $z^{\alpha_0}_{\alpha+t}y_s$ depend uniformly continuously on $\alpha$ and $s\in(t-\ve,t+\ve)$. In particular, if $s$ is sufficiently close to $t$, then both $z^{\alpha_0}_{\alpha+s}y_t$ and $z^{\alpha_0}_{\alpha+t}y_s$ are close to their common limit $z^{\alpha_0}_{\alpha+t}y_t$ uniformly in $\alpha$. It follows that the middle integral goes to $0$ for $s\to t$. Sending also $\ve\to0$, we see that the left dilation is continuous.\qed

\lf
Proposition \ref{cdprop} is more than what we actually need for continuity of the \nbd{E_0}semigroup, but it shows that we have also a reasonable notion of continuous left dilation.

\bcor
The \nbd{E_0}semigroup $\vt^v$ is strictly continuous.
\ecor

\proof
This is a more elaborate version of the proof of \cite[Theorem 10.2]{BhSk00} and a couple of similar results about continuity of \nbd{E_0}semigroups we contructed in \cite{Ske01p,Ske03b,BBLS04}. We must show that for every $a\in\sB^a(E)$ and every $x\in E$ the function $t\mapsto\vt^v_t(a)x$ is continuous. As usual with semigroups, it is sufficient to show continuity at $t=0$. Let $\zeta$ be the continuous section of unit vectors from Lemma \ref{cflem} and recall that $\zeta_0=\U$. In particular, for every $x\in E$ by Proposition \ref{cdprop} $v_\ve(x\odot\zeta_\ve)=x\zeta_\ve$ converges to $x\U=x$ for $\ve\to0$. Thus, taking also into account that $v_t(a\odot\id_t)=\vt^v_t(a)v_t$, we find that
\bmun{
ax-\vt^v_\ve(a)x
~=~
(ax-v_\ve(ax\odot\zeta_\ve))+(\vt^v_\ve(a)v_\ve(x\odot\zeta_\ve)-\vt^v_\ve(a)x)
\\
~=~
(ax-(ax)\zeta_\ve)+\vt^v_\ve(a)(x\zeta_\ve-x)
}\emun
is small for $\ve$ sufficiently small.
\qed

\lf
By Corollay \ref{univeccor} $E$ has a unit vector $\xi$. The only thing that remains to be shown is that the continuous structure induced by $\vt^v$ and $\xi$ is the same as the original one.

\bprop
A section $x$ is in $CS_i(E^\odot)$, if and only if $t\mapsto\vt^v_t(\xi\xi^*)x_t=\xi x_t$ is continuous.
\eprop

\proof
The forward implication is clear from Proposition \ref{cdprop}. For the backward implication we conclude indirectly. If $x$ is not locally uniformly bounded, then neither is $t\mapsto\xi x_t$, thus, this function cannot be continuous. So, we may assume that $x$ is locally uniformly bounded. Let us calculate $\abs{\xi x_t-\xi x_s}^2$ as in \eqref{secnd}. We find
\beqn{
\abs{\xi x_t-\xi x_s}^2
~=~
\int_d^{d+1}\abs{\xi_{\alpha+s}x_t-\xi_{\alpha+t}x_s}^2\,d\alpha,
}\eeqn
now for arbitrary $d\ge0$ because for $\alpha_0$ we may choose $0$. If $\ve$ is small, then $\xi_{\alpha+t}$ is close to $\xi_{\alpha+t-\delta}\xi_\delta$ uniformly in $\alpha$ and $\delta\in\SB{0,\ve}$ and locally uniformly in $t$. So, for $s\in\SB{t-\ve,t}$ the integral is close to
\beqn{
\int_d^{d+1}\abs{\xi_{\alpha+s}x_t-\xi_{\alpha+s}\xi_\ve x_s}^2\,d\alpha
~=~
\int_d^{d+1}\abs{x_t-\xi_\ve x_s}^2\,d\alpha
~=~
\abs{x_t-\xi_\ve x_s}^2
}\eeqn
locally uniformly in $t$. The function $\ve\mapsto i_{\ve+s}(\xi_\ve x_s)$ is continuous. So, if $\ve$ is small, then $i_{\ve+s}(\xi_\ve x_s)$ is close to $x(s)$. Thus, if $\ve$ is sufficiently small, then $\abs{\xi x_t-\xi x_s}^2$ is close to $\abs{x(t)-x(s)}$. We conlude that if $t\mapsto x(t)$ is not continuous, then neither is $t\mapsto\xi x_t$.\qed

\lf
We summarize.

\bthm\label{mthm}
Every continuous product system of correspondences over a unital \nbd{C^*}algebra is the continuous product system associated with a strictly continuous \nbd{E_0}semigroup that acts on the algebra of all adjointable operators on a unital Hilbert module.
\ethm

\section{Concluding remarks}\label{concSEC}

In these notes we discussed the simplest case of the relation between product systems and \nbd{E_0}semigroups for Hilbert modules. From every \it{strictly continuous} \nbd{E_0}semigroup acting on the algebra of all adjointable operators on a \it{unital} Hilbert module (or, more generally, on a full Hilbert module over a unital \nbd{C^*}algebra; see Remark \ref{unifullrem}) we obtain a \it{continuous} product system of correspondences over a \it{unital} \nbd{C^*}algebra and by Theorem \ref{mthm} every such product system arises in that way. A different question is in how product systems classify \nbd{E_0}semigroups.\footnote{If the \nbd{E_0}semigroups act on the same $\sB^a(E)$, then we obtain the usual classification up to cocycle conjugacy; see \cite[Theorem 2.4]{Ske02}. If two \nbd{E_0}semigroups act on two strictly isomorphic $\sB^a(E)$, so that the two Hilbert modules are \it{Morita equivalent}, then the \nbd{E_0}semigroups are cocycle conjugate, if and only if their product systems are Morita equivalent by the same Morita equivalence; see \cite[Proposition 4.7]{Ske04p}. However, it is easy to construct example of \nbd{E_0}semigroups that act on nonisomorphic $\sB^a(E)$, but have the same product systems.}

By Lemma \ref{cflem} the members in every continuous product system of correspondences over a unital \nbd{C^*}algebra have unit vectors. This allowed to adopt Arveson's point of view for the Hilbert space version in \cite{Arv06} also for modules. In all other versions we shall discuss in \cite{Ske06p1} (so far) we do not know about existence of unit vectors in the involved product systems.\footnote{Full discrete product systems, for instance, need not have unit vectors, not even in the case of von Neumann correspondences; see \cite[Examples 2.1 and 9.5]{Ske04p}.} Appart from a measurable version of the present notes (see below), these versions are: 1.) Algebraic product systems of full correspondences over a unital \nbd{C^*}algebra or of strongly full von Neumann correspondences. (Here the direct integrals will be with respect to the counting measure and, therefore, the contructed \nbd{E_0}semigroup will not be continuous.) 2.) Strongly continuous (or measurable) product systems of strongly full von Neumann correspondences. For all these versions we have to stick to the results in \cite{Ske04p} about existence of left dilations of the discrete subsystem and pass through the manipulations as indicated in \eqref{idea} in full generality.\footnote{We should emphasize that the complications in Propositions \ref{rcllprop}, \ref{denseprop} and, in particular, surjectivity of the $v_t$ in Proposition \ref{surprop} (assuring that the endomorphisms $\vt^v_t$ are unital) are not due to the construction in \cite{Arv06} but because the case of Hilbert modules is technically considerably more involved. (Also in a proof based on \eqref{idea} we have to face similar problems in showing that the three correspondences $E_t\odot\int_a^bE_\alpha\,d\alpha$ and $\int_{a+t}^{b+t}E_\alpha\,d\alpha$ and $\bfam{\int_a^bE_\alpha\,d\alpha}\odot E_t$ are isomorphic in the obvious way.) The proof in \cite{Arv06} for the analogue of Proposition \ref{surprop} for Hilbert spaces is much simpler.}

As it is our definition of continuous product system that led to unit vectors, the reader might object that this definition is too restrictive. The more important it is to underline that the relevant part of the definition is, actually, less restrictive than Arveson's. Namely, what we need in order to show existence of unit vectors in the proof of Lemma \ref{cflem} is only Property \eqref{csdef} of Definition \ref{cPSdef}. An Arveson system, appart from its structure of an algebraic product system, is a measurable bundle of Hilbert spaces $H_t,t>0$ isomorphic to the trivial bundle $(0,\infty)\times H_0$ for an infinite-dimensional separable Hilbert space $H_0$. But this bundle is also continuous, and adding a one-dimensional subspace of $H_0$ at time $t=0$ does not change this. Property \eqref{csdef} is weaker, as we do not require that the injections $i_t,t>0$ are surjective. The only difference is that Arveson requires (more or less) that products of measurable sections are measurable, while we require that products of continuous sections are continuous. In fact, the measurable version we will treat in \cite{Ske06} will simply replace the condition about continuity of products of continuous sections by measurability. Property \eqref{csdef} remains unchanged!

As far as von Neumann versions are concerned, for a normal \nbd{E_0}semigroup $\vt$ the family of projection $\vt_t(\xi\xi^*)$ will no longer be strictly continuous but only strongly continuous (in the strong topology of von Neumann modules). Consequently, in Property \eqref{csdef} continuous sections will be replaced by strongly continuous sections. But then the arguments in the proof of Lemma \ref{cflem} that led to unit vectors do no longer work. (The invertibles are not open for the strong topology.)

As far as measurability is concerned we would like to say that in \cite{Ske06} we reduced the problem to measurability of a certain unitary group. (This is much easier to treat than the continuity problem of a proper \nbd{E_0}semigroup by applying the standard result \cite[Theorem 10.2.3]{HiPhi57} of Hille and Phillips. See, for instance, the proof of \cite[Proposition 2.5(i)]{Arv89}.) This procedure, which we explain very briefly, works also for modules. The basis is to contruct not only a left dilation of the product system but also a \hl{right dilation}, that is, a Hilbert space $H$ with a faithful nondegenerate representation of $\cB$ and a family of left linear unitaries $w_t\colon E_t\odot H\rightarrow H$ iterating associatively. Then $u_t:=(v_t\odot\id_H)(\id_E\odot w_t^*)$ defines a unitary group acting on the Hilbert space $E\odot H$, which gives back $\vt$ as the restriction of the automorphism semigroup $u_t\bullet u_t^*,t\ge0$ to $\sB^a(E)\odot\id_H\subset\sB(E\odot H)$.\footnote{The fact that this semigroup is nontrivial shows that the elements $xy_t\odot z$ and $x\odot y_tz$ in $E\odot H$, in general, are different. So, the extension of the product to elements in the spaces of the left and the right dilation is no longer associative.} As for left dilations, for existence of a right dilation we need a right dilation of the discrete subsystem of $E^\odot$. For \nbd{C^*}correspondences this is the result of Hirshberg \cite{Hir05a}, while for von Neumann correspondences this is our result \cite[Theorem 7.6]{Ske04p}. In both cases the members $E_n$ of the product system must have faithful left action of $\cB$, so generality is slightly reduced.

Last but not least, we mention that, actually, Arveson \cite{Arv06} constructed a right dilation of the product system, while in \cite{Ske06} we constructed a left dilation. For Hilbert spaces there is no problem in switching from left to a right dilation simply by reversiong in all tensor products the order of the factors. (In fact, this is what we did in \cite{Ske06p4} in order to compare \cite{Arv06} with \cite{Ske06}.) Nothing like this is possible for modules! ($E\odot E_t$ has no meaning. And also the stability condition for section $x_{\alpha+1}=\zeta_1x_\alpha$, written in the reverse order $x_{\alpha+1}=x_\alpha\zeta_1$ would produce nonsense for the definition of the semiinner product on $\sS$.) This, clearly, underlines that the ``correct'' product system associated with an \nbd{E_0}semigroup is the one that is connected with the \nbd{E_0}semigroup by a left dilation, not by a right dilation. In fact, a right dilation gives rise to a \it{nondegenerate representation} $\eta_t(x_t)\colon h\mapsto x_th$ of the product system, and this is what Arveson constructed. Such a representation gives rise to an \nbd{E_0}semigroup on the von Neumann algebra $\sB^{bil}(H)$. Only in the von Neumann case, this algebra may be suitably interpreted as algebra of operators on a von Neumann module, but a von Neumann module over the commutant of $\cB$. In this case, the product system may be recovered as a family of intertwiner spaces. This is part of a far reaching duality between a von Neumann correspondence and its commutant we introduced in \cite{Ske03c}. (This relation is explained in \cite[Sections 7 and 8]{Ske04p} and in \cite{Ske06p3}.) In the \nbd{C^*}case, this \nbd{E_0}semigroup does not give back uniquely the product system for which we give a right dilation! Therefore, in the \nbd{C^*}case only the \nbd{E_0}semigroups coming by left dilations and not the \nbd{E_0}semigroup coming by nongegenerate representations, that is, by right dilations, have a ``good'' relationship to product systems.

\setlength{\baselineskip}{2.5ex}


\newcommand{\Swap}[2]{#2#1}\newcommand{\Sort}[1]{}
\providecommand{\bysame}{\leavevmode\hbox to3em{\hrulefill}\thinspace}
\providecommand{\MR}{\relax\ifhmode\unskip\space\fi MR }
\providecommand{\MRhref}[2]{%
  \href{http://www.ams.org/mathscinet-getitem?mr=#1}{#2}
}
\providecommand{\href}[2]{#2}


\end{document}